\documentclass[preprint,review,12pt]{elsarticle}
\usepackage{amssymb, amsmath}
\usepackage{epsfig}
\usepackage{subfigure}

\newtheorem{thm}{\sc Theorem}[section]																						%
\journal{}

\begin{document}

\begin{frontmatter}

\title{A remark for dynamic equations on time scales}

\author[ma]{M. Akhmet\corref{cor1}}
\ead{marat@metu.edu.tr}
\ead[url]{http://www.metu.edu.tr/~marat}
\cortext[cor1]{Corresponding author. Tel.: +90 312 210 53 55; fax: +90 312 210 29 72.}
\address[ma]{Department of Mathematics, Middle East
Technical University, 06531 Ankara, Turkey}

\author[mt]{M. Turan}
\ead{mturan@atilim.edu.tr, mehmetturan21@gmail.com}
\ead[url]{http://www.atilim.edu.tr/~mturan}
\address[mt]{Department of Mathematics, 
Atilim University, 06836, Incek, Ankara, Turkey}

\begin{abstract}
We give a proposal to generalize  the concept of the differential equations on time scales 
\cite{H88, AH88, H90}, such that they can be more appropriate for the analysis of 
real world problems, and  give more opportunities  to  increase the  theoretical depth 
of investigation.  
\end{abstract}

\begin{keyword} Differential equations on time scales; Transition condition; Modeling
\end{keyword}

\end{frontmatter}

\section{Differential equations on time scales in Hilger's sense and  the first  remark}
Let  us remind  the differential equations  on time scales, proposed by
Hilger \cite{H88, AH88, H90}. The main element of the equations is the time scale, 
which is understood as a nonempty closed subset, $\mathbb T,$ of the 
real numbers. On a time scale ${\mathbb T},$ the functions $\sigma(t):=\inf\{s\in{\mathbb T}: s>t\}$ 
and $\rho(t):=\sup\{s\in{\mathbb T}: s<t\}$ are called the forward 
and backward jump operators, respectively. The point $t\in{\mathbb T}$ 
is called right-scattered if $\sigma(t)>t,$ and right-dense if $\sigma(t)=t.$ 
Similarly, it is called left-scattered if $\rho(t)<t,$ and left-dense if $\rho(t)=t.$

The $\Delta-$derivative of a continuous function $\varphi,$ at a right-scattered 
point is defined as 
\begin{eqnarray}\label{dd1}
\varphi^\Delta(t):=\frac{\varphi(\sigma(t))-\varphi(t)}{\sigma(t)-t},
\end{eqnarray}
and at a right-dense point it is defined as
\begin{eqnarray}\label{dd2}
\varphi^\Delta(t):=\lim_{s\rightarrow t}\frac{\varphi(t)-\varphi(s)}{t-s},
\end{eqnarray}
if the limit exists.

A differential equation 
\begin{equation}\label{tsde}
 y^\Delta(t)=f(t,y), \quad t\in{\mathbb T}
\end{equation}
is said to be a differential equation on time scale, where function 
$f(t, y):{\mathbb T}\times{\mathbb R}^n\rightarrow {\mathbb R}^n$ in 
(\ref {tsde}) is assumed to be rd-continuous on ${\mathbb T}\times{\mathbb R}^n$ \cite{BP}.

\begin{figure}[htb]  
	\centering 
	\includegraphics[height=50mm]{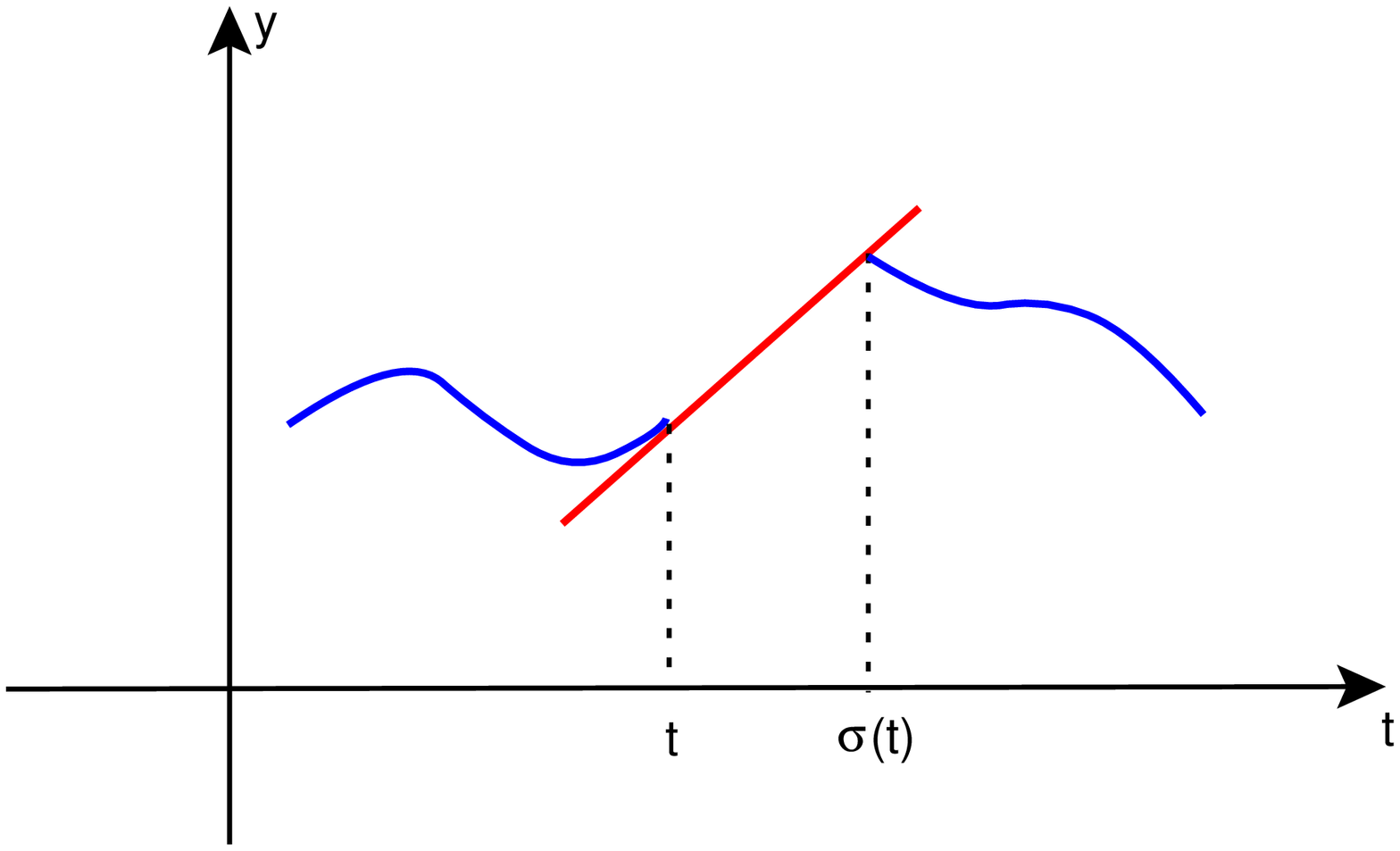}
	\caption{}\label{fig1}
\end{figure}

\begin{figure}[htb]  
	\centering 
	\includegraphics[height=50mm]{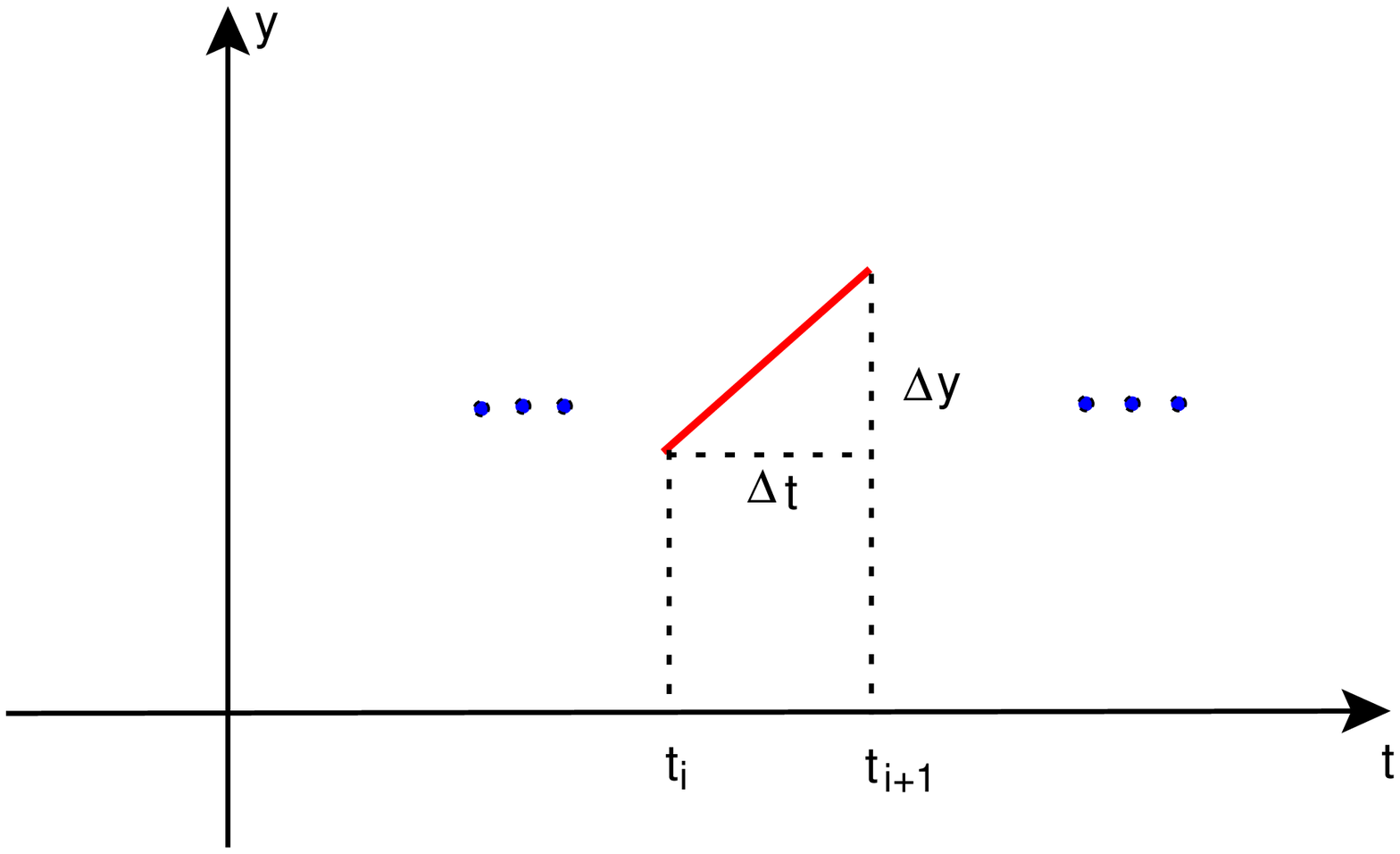}
	\caption{}\label{fig2}
\end{figure}

If a point $t$ is right-dense, then the $\Delta$-derivative is the ordinary 
derivative. Otherwise, $t\in{\mathbb T}$ is a right-scattered point, 
and the $\Delta$-derivative is defined by means of the quotient.
Let us discuss the second case more attentively. By the definition 
of the $\Delta-$derivative, assuming that $t$ is not an isolated point, 
one has that  $\varphi(\sigma(t)) =\varphi(t)+\varphi'(t)(\sigma(t) - t).$  
Now, compare this expression with the numerical approximation of the 
function on the interval $[t, \sigma(t)].$ Setting 
$\Delta t  = \sigma(t) - t$ and $\Delta \varphi = \varphi(\sigma(t))-\varphi(t),$ 
we have that $\Delta \varphi  = \varphi'(t)\Delta t.$  Thus, one can see  that  
the idea of the $\Delta-$derivative, as well as of the $\nabla-$derivative is  
close to  the basic one  for the numerical analysis \cite{SB}. 
See Figures \ref{fig1} and \ref{fig2} for illustration. 

That is, in its heart the concept of the differential equations 
on time scales is closely adjoint with numerical methods. However, 
for numerical methods, of course as far as the definition in 
$(\ref{dd1})$ is considered, smaller $\sigma(t)-t$ gives a better 
approximation. In the former theory of differential equations 
on time scales, the request was not considered, and the elements 
of linear change over the intervals near right-scattered 
(left-scattered) points are strong there. To return the 
non-linearity sense in the discussion, one needs to modify the 
equations. In our papers \cite{A05a, AT06, AT09a}, the concept 
of transition conditions were proposed, when we consider the 
simple version of the time scale as a union of disjoint sections. 
In \cite{A05a}, these systems emerge as an auxiliary instrument 
of a discontinuous dynamical system investigation. In Section 
\ref{sec3}, we generalize this idea for all possible time scales.

\section{Modeling proposals}\label{sec2}
In this part we shall answer two questions: Why the differential 
equations on time scales are needed by real world problems? 
How the Hilger's equations can be extended for these? Consider 
a real world process, which can be symbolically given as a function 
$\psi(t),$ where $t$ is from a time scale. Suppose for this moment 
that the time scale admits only non-isolated moments. If the time 
is increasing, we have to prescribe what the value at $\sigma(t)$ 
is equal to for the right-scattered point $t.$ In Hilger's theory, 
the value is evaluated by using the derivative of $\psi$ at the moment 
$t,$ and this is very restrictive, since the difference $\sigma(t)-t$ 
is not supposed to be small! It is reasonable to assume that there 
exists more general law of transition for $\psi$ from the moment 
$t$ to $\sigma(t).$ In general, it can be given with a function 
$\psi(\sigma(t)) = J(t, \psi(t)).$ We call this expression as the 
transition condition. To illustrate, consider a population developing 
in time. We may assume that the size satisfies a local law in some 
borders, say in summer time, and it is not under the law in winter time. 
Moreover, the local law in winter time is not known. Nonetheless, we 
can evaluate the change in the time, that is the size at the beginning 
of the season and that at the end of the season. Assume that the change 
satisfies a new formula, which is naturally a discrete one. In that case, 
we come to the idea of the equations on time scales, but more general 
than those investigated before.

\section{Differential equations on time-scales with  transition conditions}\label{sec3}

Consider  the time scale $\mathbb T$ given above, and the same  $\Delta-$derivative of a continuous function $\varphi,$ defined  by  (\ref{dd1}) and (\ref{dd2}).  

 Denote  by ${\mathbb T}_{\sigma}$ the  set of all right-scattered points of ${\mathbb T}.$
Consider a function $f(t, y):{\mathbb T}\backslash {\mathbb T}_{\sigma}\times{\mathbb R}^n\rightarrow {\mathbb R}^n,$  which  is assumed to be rd-continuous on ${\mathbb T}\backslash {\mathbb T}_{\sigma}\times{\mathbb R}^n.$ Define a function $J(t, y):{\mathbb T}_{\sigma}\times{\mathbb R}^n\rightarrow {\mathbb R}^n.$   Now, on the time scale, we introduce a piecewise function $F(t,y),$ which  is equal to  the function  $f(t, y),$ if $t$ is a right-dense point, and it  is equal to $J(t, y),$ if $t$ is a right-scattered  point. We suggest  to  consider the following  differential equation on the time scale 
\begin{equation}\label{tsdenew}
 y^\Delta(t)=F(t,y), \quad t\in{\mathbb T}.
\end{equation}

Condition that  all solutions are  continuous on the time scale  should be assumed. 

To demonstrate  the opportunities, which  one can have   with the newly  introduced  equations, consider the following theorem of existence and uniqueness.
\begin{thm} Let ${\mathbb T}$ be a time scale, $t_0\in{\mathbb T},$ $y_0\in{\mathbb R}^n,$ $a>0$ with $\inf{\mathbb T}\leq t_0-a$ and $\sup {\mathbb T}\geq t_0+a,$ and put
$$I_a=(t_0-a, t_0+a) \quad and \quad U_b=\{y\in{\mathbb R}^n: |y-y_0| < b\}.$$
Suppose that $f:I_a\times U_b \to {\mathbb R}^n$ is rd-continuous, bounded with bound $M>0,$ Lipschitz continuous with Lipschitz constant $L>0,$ and $$\left\|\frac{J(t,y)}{\sigma(t)-t}\right\| \leq N$$ for all  $t\in{\mathbb T}_{\sigma}\cap I_a$ and $y\in U_b.$ Then the IVP
\begin{eqnarray}\label{eqnEU}
y^{\Delta}=F(t,y), \quad y(t_0)=y_0,
\end{eqnarray}
has exactly one solution on $[t_0-\alpha, t_0+\alpha],$ where
$$\alpha=\min\left\{a, \frac{b}{\max\{M, N\}}, \frac{1-\varepsilon}{L}\right\} \quad for \  some \ \varepsilon >0.$$
If $t_0\in{\mathbb T}_\sigma$ and $\alpha < \sigma(t_0)-t_0,$ then the unique solution exists on the interval $[t_0-\alpha, \sigma(t_0)].$ \label{thmEU}
\end{thm}

One  could compare the  last  assertion with  the following, from \cite{BP}.
The condition on the transition function $J$ is added, in this particular case. In general,  properties of the function  will  generate new theoretical, and what  is most  important,  application prospects for  these equations.  One can  see that  many  results in the theory, which  have been obtained earlier can, now, be generalized,  for newly  proposed systems with  additional interesting properties that transition conditions can  admit.   

Another opportunity   is to  consider state-dependent  time scales \cite{AT09a}.   Let   a set $\mathcal D$ in $\mathbb R \times \mathbb R^n,$  be given, such  that  for each  $x \in {\mathbb R}^n,$  the projection of the intersection,  ${\cal D}_x,$ of the line $(t,x), t  \in \mathbb R$ with $\cal D$ on the time axis  is  a time scale in the sense of Hilger. Thus, one can define  the function $\sigma(t,x) = \sup\{s\in {\mathcal D}_x: s > t\}$  for each  point $(t,x) \in {\cal D}.$  We assume that  the function $\sigma$  is continuous in $x.$ If  the moment $t$ is such  that  $(t,\varphi(t))$ is right-dense in the correspond time scale, ${\cal D}_{\varphi(t)},$ then the $\Delta-$derivative  is defined by  (\ref{dd2}). Otherwise  it is equal to 
\begin{eqnarray}\label{dd1a1}
\varphi^\Delta(t):=\frac{\varphi(\sigma(t,\varphi(t)))-\varphi(t)}{\sigma(t,\varphi(t))-t}.
\end{eqnarray}
Now, define a function $G(t,x),$ which  is equal to  a continuous function $f(t,x)$ at a right-dense point $(t,x),$ and to function $J(t,x),$ if the point is a right-scattered one.  Then, one can  discuss the following equation 
\begin{equation}\label{tsdenewspace}
 y^\Delta(t)=G(t,y), \quad (t,y) \in {\cal D}.
\end{equation}
It  is obvious that  the equations  with state-dependent  time scales can  also  play a certain role in this theory,  and they  will be used in the modeling of real world problems. The analysis may  be of great interest, since new opportunities related to  the non-linear properties  of the time scales  will emerge.  

\section{Unification of differential  equations with  discontinuities}
Finally, we want to mention that   very  profitable ways  for  differential  equations on time scales is to   investigate them together with  different  types of discontinuities. These days there are several well  extended theories such  that  systems with  discontinuities in the right-hand-side, arguments, impulsive differential equations. They  are realized in very  interesting application projects as well \cite{awr}-\cite{w1}. We are sure that  integration of the methods and results of these theories  with  above mentioned extensions of differential equations on time scales can  provide really  interesting horizons,  which  can  give a light on new ways  for  modeling    processes in mechanics, electronics, medicine,  biology. One of the particular results in this sense is our recent paper \cite{AT06}, where certain class of  differential equations on time scales  is embedded in  differential equations with  fixed moments of impulses.  

\end{document}